\documentclass[11pt]{amsart}
\usepackage[margin=1in]{geometry}
\usepackage{graphicx}
\usepackage[dvipsnames]{xcolor}
\usepackage{amssymb}
\usepackage{enumitem}
\usepackage{array,multirow}
\usepackage[all]{xy}
\usepackage{amsmath}
\usepackage{amsthm}
\usepackage{hyperref}
\usepackage[nameinlink]{cleveref}

\theoremstyle{plain}

\newtheorem{thm}{Theorem}[section]
\newtheorem{cor}[thm]{Corollary}
\newtheorem{prop}[thm]{Proposition}

\theoremstyle{definition}
\newtheorem{quest}[thm]{Question}

\newtheorem{remark}[thm]{Remark}
\newtheorem{exl}[thm]{Example}

\hypersetup{colorlinks=true,linkcolor=Cyan,citecolor=Violet}

\newcommand{\N}{{\ensuremath{\mathbb{N}}}}
\newcommand{\Z}{{\ensuremath{\mathbb{Z}}}}
\newcommand{\Q}{{\ensuremath{\mathbb{Q}}}}

\newcolumntype{C}{>{\boldmath\textbf}c}

\begin{document}
\title{Amphichiral knots with large 4-genus}

\author[A.~N.~Miller]{Allison N. Miller}
\address{Department of Mathematics, Rice University, Houston, Texas, USA}
\email{allison.miller@rice.edu}

\def\subjclassname{\textup{2020} Mathematics Subject Classification}
\expandafter\let\csname subjclassname@1991\endcsname=\subjclassname

\begin{abstract}
For each $g>0$ we give infinitely many knots that are strongly negative amphichiral, hence rationally slice and representing 2-torsion in the smooth concordance group, yet which do not bound any locally flatly embedded surface in the 4-ball with genus less than or equal to $g$. Our examples also allow us to answer a question about the 4-dimensional clasp number of knots. 
\end{abstract}

\maketitle

\section{Introduction}

A knot $K$ in $S^3$ is called \emph{strongly negative amphichiral} if there exists an orientation reversing involution $\varphi \colon S^3 \to S^3$ such that $\varphi(K)= K$. Many concordance invariants vanish on such knots, including the classical Tristram-Levine signature function~\cite{LevineAlgSlice, Tristram} and more modern invariants coming from Heegaard Floer and Khovanov homology like the $\tau$-invariant~\cite{OS03g4},  $\nu^{+}$-invariant~\cite{HomWu16}, $\Upsilon$-invariant~\cite{OSS},   $s$-invariant~\cite{Rasmussen}, $s_n$-invariants~\cite{LobbKR, WuKR},   $s^\#$-invariant~\cite{KronheimerMrowkaS}, and  $\gimel$-invariant~\cite{LL19}. Notably, this list contains almost all known lower bounds on the \emph{4-genus}, or minimal genus of a (smoothly or locally flatly) embedded orientable surface in $B^4$ with boundary the given knot. However, we use Gilmer's bound on the topological 4-genus~\cite{GilmerGenus} coming from Casson-Gordon signatures~\cite{CassonGordon2} to prove the following. 

\begin{thm}\label{thm:main}
For any $g >0$, there exists a knot $K$ with the following properties:
\begin{enumerate}
\item $K$  is strongly negative amphichiral.
\item $K$ can be transformed to a smoothly slice knot by either (a) changing some crossings $(+)$ to $(-)$ or (b) changing some crossings $(-)$ to $(+)$. 
\item the topological 4-genus of $K$is strictly larger than $g$. 
\end{enumerate}
\end{thm}
In fact, something more is true, and proven in Proposition~\ref{prop:thework}: for any $g \in \N$ there exists an infinite family of knots $\{K^k\}_{k \in \mathbb{N}}$, generating a  subgroup of the concordance group isomorphic to $\left(\Z_2\right)^{\infty}$, such that any nontrivial sum $K= \#_{j=1}^m K^{k_j}$ satisfies the conclusions of Theorem~\ref{thm:main}.  Moreover, each of the knots $K^k$ is algebraically slice, so we incidentally reprove a result of Livingston~\cite{LivingstonOrder2} that there is a $\left(\Z_2\right)^{\infty}$-subgroup of the concordance group consisting of algebraically slice knots.\vspace{.2cm}

Negative amphichiral knots, if not slice, represent 2-torsion elements of the smooth concordance group;  a still-open question of Gordon asks whether all 2-torsion elements have such representatives~\cite[Problem 16]{ProcSemPsB}. We therefore obtain the following corollary to Theorem~\ref{thm:main}, which appears to be previously unknown.

\begin{cor}\label{cor:2tor}
There exist 2-torsion knots with arbitrarily large 4-genera. 
\end{cor}

A knot $K$ is called \emph{rationally slice} if there exists a smooth 4-manifold $W$ with boundary $\partial W=S^3$ and $H_*(W; \Q)= H_*(B^4; \Q)$ such that $K$ bounds a smoothly embedded null-homologous disc in $W$. Every strongly negative amphichiral knot is rationally slice~\cite{KawauchiRational}, and so Theorem~\ref{thm:main} also answers a question of~\cite{HKPS20} in the affirmative. 
\begin{cor}\label{cor:ratslice}
There exist rationally slice knots with arbitrarily large 4-genera. 
\end{cor}

The \emph{4-dimensional clasp number} $c_4(K)$  of a knot $K$ is the
minimal number of transverse double points across all immersions of $D^2$ in $B^4$ with $\partial D^2= K$. Similarly, $c_4^{+}(K)$ (respectively $c_4^{-}(K)$) is defined to be the minimal number of positive (resp. negative) transverse double points across all immersions of $D^2$ in $B^4$ with $\partial D^2= K$. It follows immediately from the definitions that $c_4^+ + c_4^- \leq c_4$; the figure-eight knot $4_1$ is the prototypical example of when this inequality is strict, since $c_4^+(4_1)= c_4^-(4_1)=0$ and yet $c_4(4_1)=1$. 
We answer a question of~\cite{JZclasp} by giving the first examples of knots for which $c_4(K)$ is arbitrarily larger than $c_4^+(K)+ c_4^-(K)$.

\begin{cor}\label{cor:c4}
The difference between $c_4(K)$ and $c_4^+(K) + c_4^-(K)$ can be arbitrarily large. 
\end{cor}
\begin{proof}
For $g \in \N$, let $K_g$ be a knot satisfying the conclusions of Theorem~\ref{thm:main}. By $(2)$, we have that  $c_4^+(K_g) + c_4^-(K_g)=0+0=0$, and by $(3)$ we have that 
\[ g< g_4(K_g) \leq g_4 ^s(K_g) \leq c_4(K_g), \]
noting that standard arguments show that for any knot $K$ the smooth 4-genus $g_4^s(K)$ is bounded above by $c_4(K)$. 
\end{proof}

Since Casson-Gordon signatures provide bounds on the topological 4-genus, it remains open whether one can find examples for the smooth analogue of Theorem~\ref{thm:main} as follows. 

\begin{quest}
For $g \in \N$, does there exist a topologically slice knot $K$ such that $g_4^s(K)>g$ and
\begin{enumerate}
\item $K$ is order 2 in the smooth concordance group?
\item $K$ is smoothly rationally slice?
\item $c_4^{+}(K)= c_4^{-}(K)=0$?
\end{enumerate}
\end{quest}

Recent work of Hom-Kang-Park-Stoffregen~\cite{HKPS20} has shown that $\{C_{2n+1,1}(4_1)\}_{n \in \mathbb{N}}$ generates a $\Z^{\infty}$-subgroup of rationally slice knots in the smooth concordance group. By work of~\cite{Millertopg4}, the topological 4-genus of $C_{2n+1,1}(4_1)$ equals 1 for all $n \in \mathbb{N}$, but it remains open whether the smooth 4-genus of $C_{2n+1,1}(4_1)$ is large. 
Since $2n+1$ is relatively prime to 2,  one can combine the work of this paper with the formulas for Casson-Gordon signatures of satellite knots given in~\cite{Litherland} and conclude that for our choice of $K_g$ satisfying the conclusions of Theorem~\ref{thm:main}, we have that $g_4(C_{2n+1,1}(K_g))>g$ for all $n \in \mathbb{N}$. We therefore state the following as an interesting open problem in either the smooth or topological categories. 

\begin{quest}
For any $g \in \mathbb{N}$, let $K_g$ be one of the knots given in Section~\ref{section:proof} that satisfies the conclusions of Theorem~\ref{thm:main}. For some or all $n \in \mathbb{N}$,
determine whether $C_{2n+1,1}(K_g)$ is infinite order in the concordance group. 
\end{quest}

We note that it remains open even whether $C_{2n,1}(K)$ must always be slice for strongly negative amphichiral $K$, though it is known that many such knots are not ribbon~\cite{Miyazaki}.

\begin{remark}
The key feature of Casson-Gordon signatures that allows us to use Gilmer's bound to establish Theorem~\ref{thm:main} when all other lower bounds on the 4-genus fail might initially seem like a flaw: no single signature gives a 4-genus bound or even a sliceness obstruction. 
While we avoid stating the precise definition of these invariants, we remind the reader that  $\sigma(K, \chi)\in \Q$ depends on not just the knot $K$ but a choice of map $\chi$ from the first homology of the double branched cover of $K$  to a cyclic group. The fact that $K$ is negative amphichiral implies that there is an involution $\iota$ on the set of such maps such that $\sigma(K, \iota(\chi))= - \sigma(K, \chi)$. As long as this involution is non-trivial, the negative amphichirality of $K$ does not force $\sigma(K, \chi)$ to vanish and there is still the potential to obtain a sliceness obstruction--and even a lower bound on the 4-genus--by considering the set of all such signatures. This could be considered as philosophically similar to the fact that Casson-Gordon signatures can obstruct knots from being concordant to their reverses~\cite{KirkLivingstonIMC}, though that result requires a careful analysis of additional structure that we are able to avoid. 
\end{remark}

\section*{Acknowledgments}
The author is indebted to Anthony Conway and JungHwan Park for thoughtful conversations and for suggesting the questions resolved by Corollaries~\ref{cor:2tor} and~\ref{cor:c4}, and to Chuck Livingston for asking about algebraic sliceness.  The author also gratefully acknowledges her partial support by NSF
DMS-1902880.

\section{Proof of Main Result} \label{section:proof}
Our examples are connected sums of certain satellites of the figure-eight knot. 
\begin{exl}
Let $J$ be a reversible knot and define $K(J)$ to be as in Figure~\ref{fig:examples}, where $\overline{J}$ denotes the mirror image of $J$, which since $J$ is reversible equals the concordance inverse $-J$. 
\begin{figure}[h!]
\includegraphics[height=3.5cm]{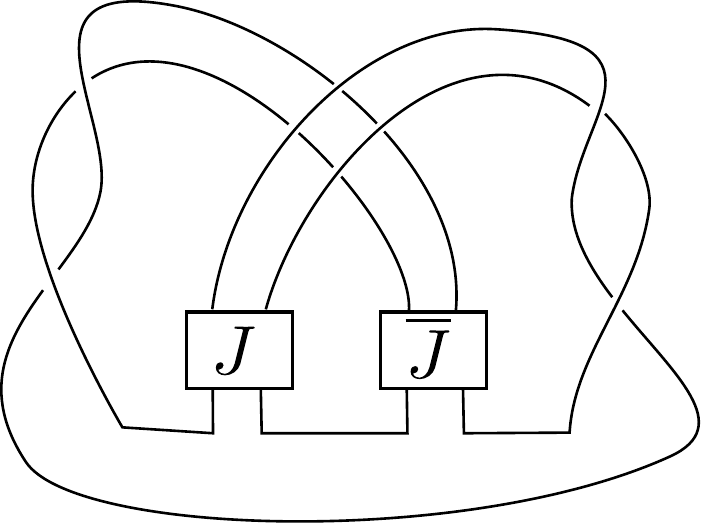} \qquad \qquad
\includegraphics[height=3.5cm]{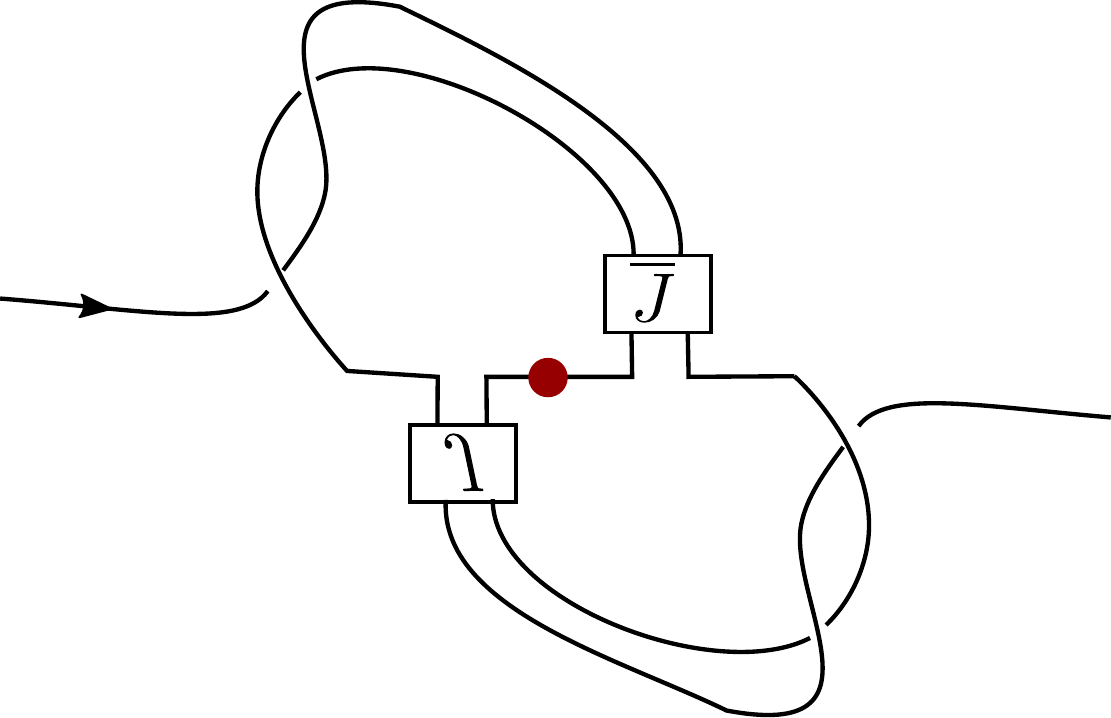}
\caption{\small The knot $K(J)$ from 2 perspectives. } \label{fig:examples}
\end{figure}
The right side of Figure~\ref{fig:examples} demonstrates that $K(J)$ is strongly negative amphichiral: rotation by 180 degrees in the plane about the marked point followed by reflection in the plane of the page takes $K(J)$ to itself. Also observe that the disc-with-bands Seifert surface for $K(J)$ visible on the left of Figure~\ref{fig:examples} demonstrates that $K(J)$ shares a Seifert form with the figure-eight knot $K_0$. \end{exl}

\begin{prop}
If $J$ is a reversible knot, then $K(J)$ has $c_4^+(K_J)= c_4^-(K_J)=0$. 
\end{prop} 
\begin{proof} Consider the knots $K_{\pm}$ as depicted in Figure~\ref{fig:crossingchange}, shown with genus one Seifert surfaces $F_{\pm}$  in disc-with-bands position. Observe that $K_+$ (respectively $K_-$) is obtained from $K_J$ by changing a single negative (resp. positive) crossing to a positive (resp. negative) crossing.  
\begin{figure}[h!]
\includegraphics[height=3.5cm]{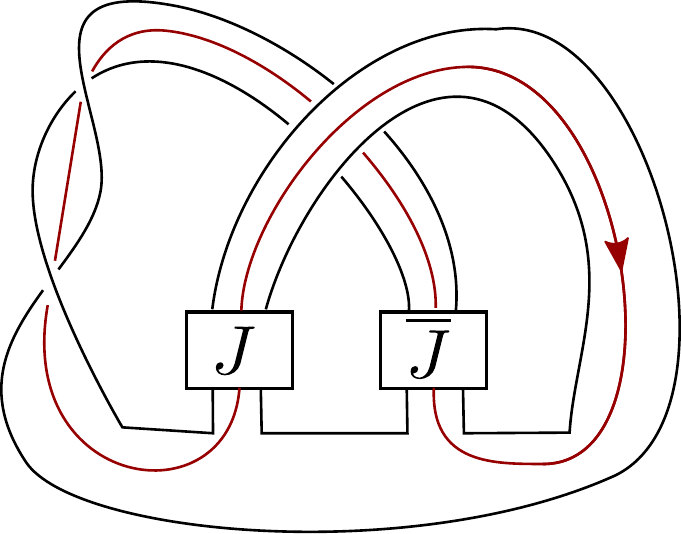}\qquad \qquad
\includegraphics[height=3.5cm]{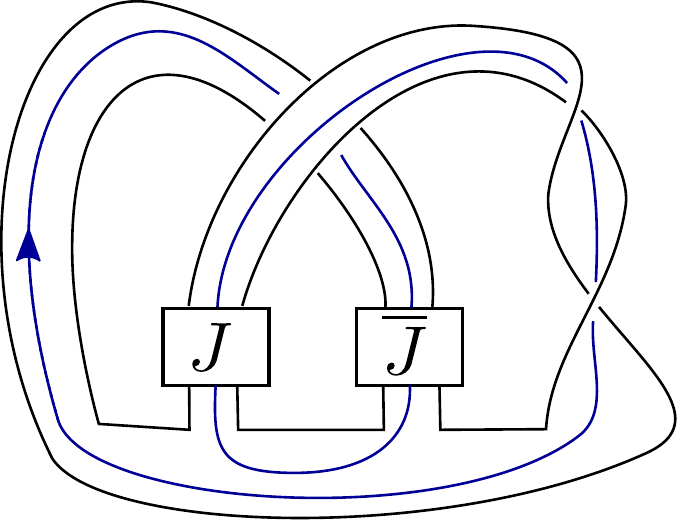}
\caption{\small $K_+$, obtained by changing a crossing from $-$ to $+$ (left) and $K_-$, obtained by changing a crossing from $+$ to $-$ (right). }\label{fig:crossingchange}
\end{figure}
Figure~\ref{fig:crossingchange} also depicts a curve $\gamma_{\pm}$ on $F_{\pm}$. Note that each of $\gamma_{\pm}$ represents a  nontrivial element of $H_1(F_{\pm})$ and is 0-framed by $F_{\pm}$; i.e. is an \emph{derivative curve}. Considered as a knot, $\gamma_+$ is $J \# \overline{J}$; since $J$ is reversible this is isotopic to $J \# -J$ and hence is slice. Similarly, the knot type of $\gamma_-$ is the slice knot $J \#-J$. Therefore, surgering the Seifert surface $F_{\pm}$ along the derivative curve $\gamma_{\pm}$ yields a smooth slice disc for $K_{\pm}$. 
We can convert this single crossing change from $K(J)$ to $K_{\pm}$ into an immersed annulus in $S^3 \times I$ from $K(J)$ to $K_\pm$. Capping each of these annuli with a smooth slice disc for $K_{\pm}$ yields the desired immersed discs bounded by $K(J)$, each with a single singularity  of different sign. 
\end{proof}

\subsection{Background results}

 For $n \in \mathbb{N}$ and a knot $K$, we let $\Sigma_n(K)$ denote the $n$th cyclic branched  cover of $S^3$ along $K$. To a knot $K$ and a map $\chi \colon H_1(\Sigma_n(K)) \to \Z_q$ one can associate the Casson-Gordon signature $\sigma(K, \chi) \in \Q$~\cite{CassonGordon2}. 
 We avoid giving the technical definition of these invariants, noting only that  they are defined in terms of the twisted intersection form of some 4-manifold  and are notoriously difficult to compute precisely. We remark for those familiar with Casson-Gordon signatures that  in the literature what we call $\sigma(K, \chi)$ is just $\sigma_1 \tau (K, \chi)$ instead. 

Our lower bound on the topological 4-genus of a knot comes from the following result of Gilmer. 
\begin{thm}[\cite{GilmerGenus}]\label{thm:gilmer}
Suppose that $K$ is a knot with $g_4(K) \leq g$. Then there is a decomposition $H_1(\Sigma_2(K))= A_1 \oplus A_2$ such that:
\begin{enumerate}
\item$ A_1$ has a presentation with at most $2g$ generators.
\item There is some $B \leq A_2$ with $|B|^2= |A_2|$ 
 such that for any prime power order $\chi \colon H_1(\Sigma(K)) \to \Z_q$, we have 
\[ |\sigma(K, \chi) + \sigma(K)| \leq 4g.\]
\end{enumerate}
\end{thm}
We remark for later that in our applications of Theorem~\ref{thm:gilmer} $K$ will always be negative amphichiral and hence have $\sigma(K)=0$. \\

Litherland proved a much more general formula for the Casson-Gordon invariants of satellite knots, but we will only need the following special case.

\begin{thm}[\cite{Litherland}]\label{thm:litherland}
Suppose $P$ is a pattern of winding number 0 described by an unknot $\eta$ in the complement of $P(U)$.  Let $x$ denote the homology class of one of the lifts of $\eta$ to $\Sigma_2(P(U))$.
For any knot $J$, there is an isomorphism $\alpha \colon H_1(\Sigma_2(P(J))) \to H_1(\Sigma_2(P(U)))$ such that for any $\chi \colon H_1(\Sigma_2(P(U))) \to \Z_q$ we have 
\[
\sigma(P(J), \chi \circ \alpha)= \sigma(P(U), \chi)+ 2\sigma_J(\omega_q^{\chi(x)}),
\]
where  $\omega_q=e^{2 \pi i/q}$ and $\sigma_J$ denotes the Tristram-Levine signature function. 
\end{thm}

As well as the knot invariant $\sigma(K, \chi)$, Casson-Gordon introduced a signature invariant $\sigma(M, \phi)$ associated to a 3-manifold $M$ and a character $\phi \colon H_1(M) \to \Z_q$. 
We will need a formula due to Cimasoni-Florens for the Casson-Gordon signature of a 3-manifold in terms of the colored signature function of a surgery link. Although this result is proved in much more generality, we state it only for the case of interest: when $M$ is obtained by surgery on a Hopf link. We thereby avoid going into the technical details of the definition of the colored signature function, noting only for the experts that the cell complex consisting of 2 discs meeting in a single arc and bounded by the Hopf link  is a C-complex in the sense of~\cite{CimasoniFlorens}, and the contractibility of this complex immediately implies that the colored signature function of the Hopf link is identically zero.

\begin{thm}\cite[Theorem 6.7]{CimasoniFlorens} \label{thm:cf}
Suppose that a 3-manifold $M$ is obtained by surgery on a Hopf link $L$ with linking matrix $\Lambda= \left[ \begin{array}{cc} a & 1 \\1 & b \end{array} \right]$. Let $q$ be prime and $\chi \colon H_1(M) \to \Z_q$ be a character such that the two meridians $\mu_{1}, \mu_2$  of $L$ are sent to nonzero elements of $\Z_q$. For $i=1, 2$ let $n_i \in \{1, \dots, q-1\}$ be the unique value satisfying $n_i \equiv \chi(\mu_i) \mod q$. 
Then
\[ \sigma(M, \chi)= -1 - \text{sign}(\Lambda) + \frac{2}{q^2}
  \left[ \begin{array}{c} n_1 \\ n_2 \end{array} \right]^T\cdot
    \left[ \begin{array}{cc} a & 1 \\ 1 & b \end{array} \right]\cdot
      \left[ \begin{array}{c} q-n_1 \\ q-n_2\end{array} \right]
\]
\end{thm} 

\subsection{Proof of Theorem~\ref{thm:main}}
We now apply Theorems~\ref{thm:litherland} and~\ref{thm:cf} to obtain a formula for the Casson-Gordon signatures of $K_J$ in terms of the Tristram-Levine signatures of $J$. 

\begin{exl}\label{exl:cgcomp}
Let $K_0$ denote the figure-eight knot. 
Note that $K(J)$ is obtained from $K_0$ by two infections along curves $\eta_1$ and $\eta_2$, as depicted in Figure~\ref{fig:infection}.
\begin{figure}[h!]
 \includegraphics[height=3cm]{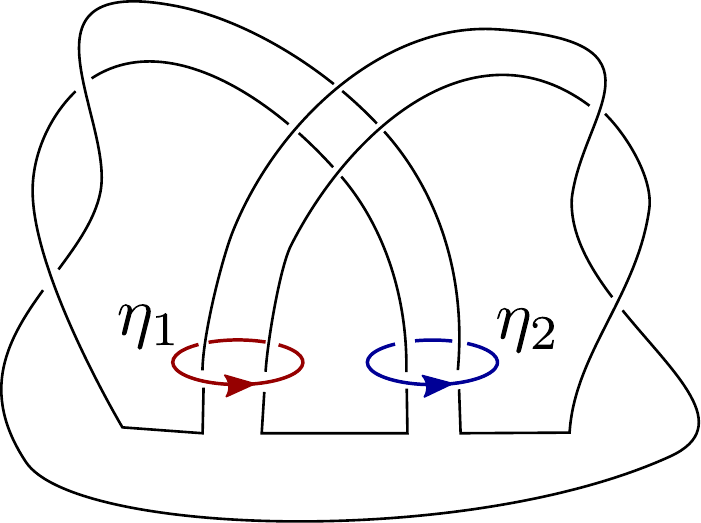}
\caption{\small The knot $K(J)$ is an iterated satellite of the figure-eight knot. } \label{fig:infection}
\end{figure}

By twice applying Theorem~\ref{thm:litherland}, we see that for any knot $J$  there is an isomorphism $\alpha \colon H_1(\Sigma_2(K(J))) \to H_1(\Sigma_2(K_0))$ such that for any character $\chi \colon H_1(\Sigma_2(K_0)) \to \Z_q$ we have 
\begin{align*}
\sigma(K(J), \alpha \circ \chi) &= \sigma(K_0, \chi)+ 2 \sigma_J(\omega_q^{\chi(\widetilde{\eta_1})})+2 \sigma_{\overline{J}}(\omega_q^{\chi(\widetilde{\eta_2})}) = \sigma(K_0, \chi)+ 2 \sigma_J(\omega_q^{\chi(\widetilde{\eta_1})})-2 \sigma_{J}(\omega_q^{\chi(\widetilde{\eta_2})})
\end{align*}

Since both $\eta_i$ curves are disjoint from the usual genus one Seifert surface for $K_0$, we can apply Akbulut-Kirby's algorithm of~\cite{AkbulutKirbyBranched} to obtain the following surgery diagram for $\Sigma_2(K_0)$, with lifts of $\eta_1$ and $\eta_2$ as indicated. (Note that we have only depicted one lift of each curve, since that is all we need to apply Theorem~\ref{thm:cf}.) 
\begin{figure}[h!]
 \includegraphics[height=2cm]{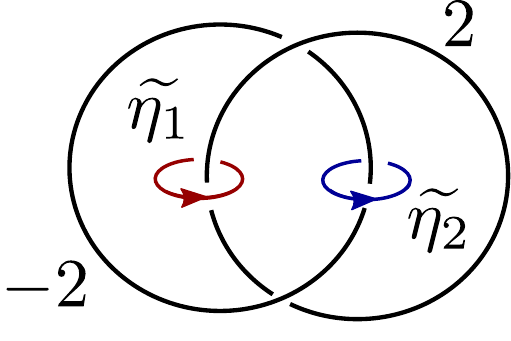}
\caption{\small A surgery diagram $L$ for $\Sigma_2(K_0)$.} \label{fig:dbc}
\end{figure}
The first homology of $\Sigma_2(K_0)$ is generated by the meridians of the components of $L$, which are isotopic to $\widetilde{\eta_1}$ and $\widetilde{\eta_2}$. The relations are given by the rows of the linking-framing matrix, and are
\[ -2[\widetilde{\eta_2}]+ [\widetilde{\eta_1}]=0 \text{ and } [\widetilde{\eta_2}] + 2 [\widetilde{\eta_1}]=0.\]
Some quick simplifications give us that $H_1(\Sigma_2(K_0)) \cong \Z_5$, generated by $a:=[\widetilde{\eta_2}]$ and such that $[\widetilde{\eta_1}] = 2[\widetilde{\eta_2}]$. 
Therefore, for any character $\chi \colon H_1(\Sigma_2(K_0)) \to \Z_5$ we have that 
\begin{align}\label{eqn:sigkj}
 \sigma(K_J, \chi \circ \alpha)= \sigma(K_0, \chi) + \sigma_J(\omega_5^{2 \chi(a)})- \sigma_J(\omega_5^{\chi(a)}).
 \end{align}

We can also use the surgery diagram of Figure~\ref{fig:dbc} to bound $|\sigma(K_0, \chi)|$.  For $j \in \Z_5$, define $\chi_j \colon H_1(\Sigma_2(K_0)) \to \Z_5$ to be the map with $\chi_j(x)=j$. Observe that $\chi_1([\widetilde{\eta_1}])=2$ and $\chi_2([\widetilde{\eta_1}])=4$. 
Therefore, Theorem~\ref{thm:cf} gives us that 
\begin{align*}
 \sigma(\Sigma_2(K_0), \chi_1)&= -1-0 +\frac{2}{25}
 \left[ \begin{array}{cc}
 1 & 2 
 \end{array} \right]
  \left[ \begin{array}{cc}
-2 & 1\\ 1 & 2 
 \end{array} \right]
  \left[ \begin{array}{c}
4 \\ 3
 \end{array} \right]= -1 + \frac{30}{25} = 1/5
\end{align*}
and 
\begin{align*}
 \sigma(\Sigma_2(K_0), \chi_2)&= -1-0 +\frac{2}{25}
 \left[ \begin{array}{cc}
 2 & 4
 \end{array} \right]
  \left[ \begin{array}{cc}
-2 & 1\\ 1 & 2 
 \end{array} \right]
  \left[ \begin{array}{c}
3 \\ 1
 \end{array} \right]
 = -1 + \frac{20}{25} = -1/5.
\end{align*}
Moreover, basic properties of Casson-Gordon signatures (or reapplying Theorem~\ref{thm:cf}) imply that $ \sigma(\Sigma_2(K_0), \chi_3)= \sigma(\Sigma_2(K_0), \chi_2)$,  $ \sigma(\Sigma_2(K_0), \chi_4)= \sigma(\Sigma_2(K_0), \chi_1)$, and $ \sigma(\Sigma_2(K_0), \chi_0)=0$.  

Since $H_1(\Sigma_2(K_0)) \cong \Z_5$ is cyclic, for any character $\chi \colon H_1(\Sigma_2(K_0)) \to \Z_5$ we have by \cite[Lemma 3 and Theorem 4]{CassonGordon2} that 
\[ |\sigma(K_0, \chi)- \sigma(\Sigma_2(K_0), \chi)| \leq 1.\]
Therefore, we conclude that for any $\chi \colon H_1(\Sigma_2(K_0)) \to \Z_5$ we have $|\sigma(K_0, \chi)| <2$. \\

\end{exl}

We are now ready to prove the following and obtain Theorem~\ref{thm:main} as a consequence. 

\begin{prop}\label{prop:thework}
Fix $g \in \mathbb{N}$. 
For  $i\in \mathbb{N}$ define $J_i= \#^{m_i} T_{2,5}$, where $m_i= 2^{2i+1} g$. 
Now, for $k \in \mathbb{N}$ define $K^k:= \#_{i=1}^{2g+2}  K(J_{k(2g+2)+i})$. 
Then $S=\{K^k\}_{k\in \mathbb{N}}$ is a collection of algebraically slice knots  such that any nontrivial sum $K= \#_{j=1}^n K^{k_j}$ satisfies the conclusions of Theorem~\ref{thm:main}. 
\end{prop}

\begin{proof}
Observe that for any choice of $J$, the knot $K(J)$ shares a Seifert form with $K_0$. Therefore, each $K^k$ shares a Seifert form with the slice knot $\#_{i=1}^{2g+2} K_0$, and hence is algebraically slice. 

Now let $K=\#_{j=1}^{n} K^{k_j}$ be a nontrivial sum of elements of $S$. We can and do assume that $k_1<k_2<\dots<k_n$.  Since conditions $(1)$ and $(2)$ of Theorem~\ref{thm:main} are preserved under connected sum, it only remains to verify condition $(3)$. 

So suppose for a contradiction that $g_4(K) \leq g$ and hence that there exists a decomposition
$H_1(\Sigma_2(K)) \cong A_1 \oplus A_2$ and a subgroup $B \leq A_2$ satisfying the conclusions of Theorem~\ref{thm:gilmer}. Let 
\begin{align*}
 \beta \colon H_1(\Sigma_2(K))  &\to \bigoplus_{j=1}^n  \left( \bigoplus_{i=1}^{2g+2}H_1(\Sigma_2(K(J_{k(2g+2)+i})))  \right)  \to  \bigoplus_{j=1}^n  \left( \bigoplus_{i=1}^{2g+2}H_1(\Sigma_2(K_0))\right)
\end{align*}
denote the isomorphism (coming from Theorem~\ref{thm:litherland} together with the additivity of Casson-Gordon signatures with respect to connected sum~\cite{Litherland}) satisfying 
\begin{align*}
\sigma\left(K, \beta \circ \left( ( \chi_i^j)_{i=1}^{2g+2}\right)_{j=1}^n \right)& = \sum_{j=1}^m \sigma(K^{k_j}, ( \chi_i^j)_{i=1}^{2g+2})\\
&= \sum_{j=1}^n \left(\sum_{i=1}^{2g+2} \sigma( K(J_{k_j(2g+2)+i}), \chi_i^j) \right) \\
&= \sum_{j=1}^n \left(\sum_{i=1}^{2g+2} \sigma( K_0, \chi_i^j)+ 2\sigma_{J_{k_j(2g+2)+i}}(\omega_5^{2 \chi_i^j(a)})-
2\sigma_{J_{k_j(2g+2)+i}}(\omega_5^{\chi_i^j(a)}) \right),
\end{align*}
where in the last equality we use Equation~\ref{eqn:sigkj} of Example~\ref{exl:cgcomp}. 

Since $H_1(\Sigma_2(K)) \cong \Z_5^{m(2g+2)}$ and $A_1$ has a presentation with at most $2g$ generators, we have that $A_1$ is isomorphic to $\Z_5^j$ for some $j \leq 2g$. Therefore $A_2$ is isomorphic to $\Z_5^{n(2g+2)-j}$ and $B$  is isomorphic to $\Z_5^{n(g+1)-j/2}$.  So $A_1 \oplus B \cong \Z_5^{j'}$ for 
\[j'= n(g+1) + j/2 \leq n(g+1)+g <n(2g+2)
\]
 and there exists a nonzero character $\chi \colon H_1(\Sigma_2(K)) \to \Z_5$ that vanishes on $A_1 \oplus B$.

The rest of the proof consists of showing that $|\sigma(K,\chi)|>4g$, using only our definition of $K$ and the hypothesis that $\chi$ is not identically zero.
Let \begin{align*} \left((\chi_i^j)_{i=1}^{2g+2}\right)_{j=1}^n:= \chi \circ \beta^{-1} \colon \bigoplus_{j=1}^n  \left( \bigoplus_{i=1}^{2g+2}H_1(\Sigma_2(K_0))\right) \to \Z_5. 
\end{align*}
 Since $\chi$ is nontrivial, there exists some $j$ such that $(\chi_i^j)_{i=1}^{2g+1}$ is not identically zero. Let $j_0$ be the maximal such $j$ and $i_0$ be the maximal $i$ such that $\chi_{i}^{j_0}$ is nonzero.
Let $\ell= k_{j_0}(2g+2)+{i_0}$. The following algebraic manipulations show that $\sigma(K(J_\ell), \chi_{i_0}^{j_0})$ so dominates the other terms that could contribute to $\sigma(K, \chi)$ that we have as desired that $|\sigma(K, \chi)|>4g$.

Recalling that $J_i= \#^{m_i} T_{2,5}$, where $m_i= 2^{2i+1} g$, we have by the additivity of Tristram-Levine signatures under connected sum that $\sigma_{J_i}(\omega_5)=\sigma_{J_i}(\omega_5^4)= -2^{2i+2} g$ and $\sigma_{J_i}(\omega_5^2)= \sigma_{J_i}(\omega_5^3)= -2^{2i+3} g$ (see KnotInfo~\cite{knotinfo} for the Tristram-Levine signature function of $T(2,5)$.)
Applying Equation~\ref{eqn:sigkj}  from Example~\ref{exl:cgcomp}, we see that 
 for any $i$ and any nonzero character $\rho \colon H_1(\Sigma_2(K(J_i))) \to\Z_5$ we have that 
\begin{align} \label{eqn:sigk}
 2^{2i+3}g -2 \leq |\sigma(K(J_i), \rho)| &=|\sigma(K_0, \rho) \pm (2\sigma_{J_i}(\omega_5)-  2\sigma_{J_i}(\omega_5^2))| \leq 2^{2i+3}g+2
\end{align}
Note that here and in the rest of the proof, we suppress the identification of each $H_1(\Sigma_2(K(J_i)))$ with $H_1(\Sigma_2(K_0))$. 

Observe that the set of natural numbers 
\begin{align}\label{eqn:obs} \{ k_{j_0}(2g+2)+i: 1 \leq i \leq i_0-1\} \cup  \bigcup_{j=1}^{j_0-1} \{k_j(2g+2)+i): 1 \leq i \leq 2g+2\} 
 \end{align}
is a subset of $\{1, \dots, \ell -1\}$, recalling that $\ell= k_{j_0}(2g+2)+i_0$. 
We therefore have that
\begin{align*}
|\sigma(K, \chi)|& =\left| \sum_{j=1}^n \sum_{i=1}^{2g+2} \sigma( K(J_{k_j(2g+2)+i}), \chi_i^j) \right| \\
&=\left| \sigma( K(J_{\ell}), \chi_{i_0}^{j_0})+
\sum_{i=1}^{i_0-1} \sigma( K(J_{k_{j_0}(2g+2)+i}), \chi_i^{j_0})
+\sum_{j=1}^{j_0-1} \sum_{i=1}^{2g+2}  \sigma(K(J_{k_j(2g+2)+i}), \chi_i^j) \right|
\\
&\geq \left| \sigma( K(J_{\ell}), \chi_{i_0}^{j_0})\right|-
\sum_{i=1}^{i_0-1} \left| \sigma( K(J_{k_{j_0}(2g+2)+i}), \chi_i^{j_0})\right|
-\sum_{j=1}^{j_0-1} \sum_{i=1}^{2g+2} \left|  \sigma(K(J_{k_j(2g+2)+i}), \chi_i^j) \right|
\\
&\geq  (2^{2 \ell+3}g-2)
- \sum_{k=1}^{\ell-1} (2^{2k+3}g+2) =:(*)
\end{align*}
where in the last inequality we use our observation from Equation~\ref{eqn:obs} together with Equation~\ref{eqn:sigk}. 
Some algebraic simplification yields that 
\begin{align*}
(*)&=8g\left(2^{2 \ell} - \sum_{k=1}^{\ell-1} 2^{2k}\right) - 2 \ell = (g/3)(2^{2 \ell+3} - 32) -2 \ell. 
\end{align*}

Now, note that since $\ell>2g+2 \geq 4$ we have that $2 \ell +3 > 11$ and so certainly $2^{2 \ell+3}-32 > 2^{2 \ell+2}$. Therefore 
\begin{align*}
|\sigma(K, \chi)| \geq (*)> (g/3) 2^{2\ell+2} - 2 \ell> 2^{2 \ell} - 2 \ell.
\end{align*}
Finally, we observe that for any $x>2$ we have $2^{2x}-2x> 2x$, since letting $f(x)= 2^{2x}-4x$ we see that   $f'(x)= \ln(4)2^{2x}-4$ is positive for all $x \geq 1$ and $f(2)= 8$. Therefore 
\[ |\sigma(K, \chi)|> 2 \ell > 4g+4>4g,\]
as desired. 
\end{proof}

\begin{remark}
The examples of Proposition~\ref{prop:thework} are far from the only knots satisfying the conclusions of Theorem~\ref{thm:main}. One could vary the base knot, for example by choosing $\{a_i\}_{i \geq 0}$ to be natural numbers such that $\{4a_i^2+1\}_{i \in \N}$ consists of pairwise relatively prime numbers. (This is easily accomplished by e.g.~letting $a_0=1$ and $a_k=\prod_{i=1}^{k-1} (4a_i^2+1)$ for $k\geq1$. ) Now, let $K_i$ be the 2-bridge knot corresponding to the rational number $\frac{4a_i^2+1}{2a_i}$, noting that indeed $K_0$ is the figure-eight knot. Choose $\{p_i\}_{i \geq 0}$ to be primes dividing $4a_i^2+1$, noting that by our choice of $a_i$ we have that $p_i$ divides $4a_j^2+1$ if and only if $j=i$.
By taking connected sums of $K_{a_i}$ analogously infected with large connected sums of $T_{2, p_i}$ and $-T_{2, p_i}$, we can essentially repeat the arguments of Proposition~\ref{prop:thework} and obtain many more  linearly independent knots satisfying the conclusions of Theorem~\ref{thm:main}. 
\end{remark}

\bibliography{bibliography}

\bibliographystyle{alpha}
\end{document}